# On evaluating the measure of strong projections in infinite dimension

Miklos Ferenczi

*Budapest University of Technology and Economics, Department of Mathematics*

**Abstract**

Projections of finite dimensional sets and their measures are investigated in infinite-dimensional power measure spaces. The starting point is the known algebraic formula, expressing the $y$-projection of a finite-dimensional set $a$ as a Boolean supremum of certain finite geometrical transformations of $a$ in the infinite-dimensional power space. This Boolean supremum somewhat unusual in classical measure theory because, it is different, in general, from the usual union of sets. The paper investigates the problem whether the power measure in the infinite-dimensional measure space is continuous with respect to the forementioned Boolean supremum. If so, then this continuity leads to a simple formula for calculating the measure of the projection of $a$. It is shown that the answer concerning this continuity is affirmative for discrete measures but false for the Lebesgue measure, for example. However, it is proved that if the concept of the $y$-projection of $a$ is replaced by that of the so-called *strong $y$-projection* of $a$, then the Lebesgue measure becomes continuous in this context and the value of the corresponding real supremum is exactly the measure of the foregoing strong $y$-projection. In this paper, the tools of the classical analysis are adapted to handle measures on Boolean algebras.

*Key words:* measurable projection, measure of projection, measures on Boolean algebras
*2020 MSC:* 28A05, 28A35

## 1. Introduction

Projections of finite-dimensional sets and their measures are investigated in infinite-dimensional power measure spaces. These projections have been extensively studied but the approach applied here is unusual in certain respects. We start from the following known algebraic property of the finite-dimensional projection in the infinite-dimensional space:

$$C_y a = \sup_{n \in \omega} \bigcup_{i=1}^{n} a(y/y_i) \qquad (1)$$

where $\mathcal{D}$ is a set algebra with unit $T^\alpha$, $\alpha \geq \omega$, $T$ is an infinite set, $a \in \mathcal{D}$, the elements of $D$ are finite-dimensional, $D$ is closed under the cylindrifications $C_y$ (identified with the appropriate projection) and under the substitutions $y/y_i$, where the dimension variable $y$ is free in $a$, the sup denotes the Boolean supremum in the set algebra $\mathcal{D}$ and $y_1, y_2, \ldots y_i, \ldots$ is an *arbitrary* $\omega$-sequence of distinct dimension variables (ordinals), not free in $a$ (see Sect. 2.1, [10]). $\mathcal{D}$ can be chosen, for example, as the set algebra $\mathcal{H}$ of the cylinders in $T^\alpha$ corresponding to the finite powers of a fixed $\sigma$-algebra $\mathcal{A}$ with unit $T$, where $T$ is countable, or as the set algebra $\mathcal{G}$ of the finite unions of the finite-dimensional rectangles in $T^\alpha$, where $T = [0, 1)$.

Assume that a probability measure $\mu$ is given on $\mathcal{D}$. The problem arises whether $\mu$ is continuous w.r.t. the supremum in (1), i.e. whether the following is true:

$$\mu(C_y a) = \sup_{n \in \omega} \mu(\bigcup_{i=1}^{n} a(y/y_i)) \qquad (2)$$

If the measure $\mu$ is continuous, then (2) serves as a simple formula to evaluate the measure of the $y$-projection of $a$ in the infinite-dimensional measure space. Notice that (2) represents a *special kind* of continuity of the measure because the Boolean supremum in (1) can not be regarded as an infinite union of sets, as is common in classical measure theory. The sets $C_y a$ and $\bigcup_{i=1}^{\infty} a(y/y_i)$ are generally *different* (see Sect. 2.1).

In Theorem 1, it is proved that the answer for the problem regarding the validity of the relation in (2), is affirmative if the measure $\mu$ is composed from the finite powers of a probability measure $p$, where $p$ is defined on a discrete measure space $\langle T, \mathcal{A}, p \rangle$ with a countable $T$, $\mu$ is defined on the algebra $\mathcal{H}$ as described above. By Theorem 2, (2) is also true if $\mu$ is the finite-dimensional Lebesgue measure $\lambda$ on the set algebra $\mathcal{G}$ of the finite unions of the finite dimensions rectangles. However, as in Sect. 3 examples show, (2) may be false if $a$ is a $n$-dimensional Lebesgue-measurable set and $\mu$ is the Lebesgue measure. It is a well-known result that the finite-dimensional Lebesgue measure $\lambda$ is continuous on $\mathcal{L}^{cyl}$ in the classical sense, where $\mathcal{L}^{cyl}$ is the set algebra of the cylinders of the finite-dimensional Lebesgue sets in $[0,1)^\alpha$([2]). However, this kind of continuity does not concern the Boolean sup in (1).

In Sect. 5, the concept of *strong projection* $\widetilde{C_y}a$ of $a$ is introduced. Strong $y$-projection means that a point $r$ belongs to the strong $y$-projection of a set $a$ ($a \subset T^\alpha$) if and only if the $r$-conditional measure of the set of $a$-points whose ordinary $y$-projection is $r$ is *positive*, i.e. this conditional distribution w.r.t. $r$ is non-zero. This projection depends on measure. I is closely related to the concept of probability quantifier $\widetilde{\exists x}\, b(x)$, where the latter is true if and only if the probability of the set $\{t : b(t) \text{ is true}\}$ is positive, i.e. non-zero. In the discrete case, the classical projection and the strong coincide.

As a consequence of Theorem 4, it is shown, that if $\mathcal{D}$ is the set algebra $\mathcal{L}^{cyl}$ of the cylinders of the finite-dimensional Lebesgue sets in $[0,1)^\alpha$, then $\mathcal{D}$ is *closed* under strong projections and the following relation holds:

$$\lambda(\widetilde{C_y}a) = \sup_{n \in \omega} \lambda(\bigcup_{i=1}^{n} a(y/y_i)) \qquad (3)$$

where (3) is obtained from (2) replacing the projection $C_y$ by the strong projection $\widetilde{C_y}$ and $\lambda$ is the finite-dimensional Lebesgue measure on $\mathcal{D}$. Thus, (3) is a formula for calculating the measure of the strong projection using the finite dimensional Lebesgue measure.

The results remain valid if the Lebesgue measure space is replaced by any other continuous probability measure space.

The measures investigated here play an important role in the theory of stochastic processes (see Kolmogorov's consistency theorem, [2]), so the results can also be adapted to that context.

## 2. Concepts

*2.1. Basic operations and structures*

Let $\mathcal{F}$ be a set algebra with unit $T^\alpha$, where $\alpha$ is a fixed infinite ordinal, $T$ is an infinite set, called the *base* of $\mathcal{F}$. For $T$, often the interval $[0,1)$, or a countable set is chosen. Let $x$, $y$, $y_k$, $i$, $k$, $n$ denote ordinals. The ordinal $k$ or $y_k$ are often called also the $k$th *dimension* or the $k$th variable.

If $r \in T^\alpha$, then let $r_k$ denote the $k$th member of $r$, where $k < \alpha$. Let $r_u^k$ denote the element of $T^\alpha$ such that $r_k$ in $r$ is changed by $u$, $u \in T, k < \alpha$. Let $\mathcal{P}(T^\alpha)$ denote the power set of $T^\alpha$.

The operation $k$th *cylindrification* is defined for the sets $a \in \mathcal{P}(T^\alpha)$ is as follows:

$$C_{y_k} a = \left\{ r : r \in T^\alpha,\ r_u^k \in a \text{ for some } u \in T \right\} \qquad (4)$$

Geometrically, it means forming the $y_k$-*cylinder* of $a$ (parallel to the $y_k$-dimension). The $y_k$-cylinder of $a$ *is associated with the* $y_k$-*projection* of $a$ (parallel to the $y_k$-dimension). From the viewpoint of logic, $C_{y_k}$ can be associated with the *existential quantification* $\exists y_k$. If a set algebra $\mathcal{F}$ with unit $T^\alpha$ is closed under the operations $C_{y_k}$'s, $k < \alpha$, then $\mathcal{F}$ is called an $\alpha$-dimensional diagonal-free *cylindric set algebra* ([1], [10]).

The *dimension set* of a set $a \in \mathcal{P}(T^\alpha)$ consists of those dimension variables $y_{k_1}, y_{k_2}, \ldots y_{k_n}, \ldots$ for which $C_{y_k} a \neq a$. It is denoted by $\triangle a$. The set $a$ is sometimes denoted as $a(y_{k_1}, y_{k_2}, \ldots y_{k_n} \ldots)$, where $y_{k_1}, y_{k_2}, \ldots y_{k_n}, \ldots$ are called also the *free variables* of $a$, the other variables are the bounded variables of $a$. If $a$ is finite-dimensional with $\triangle a = \{y_{k_1}, y_{k_2},\ \ldots y_{k_n} : k_i < \alpha\}$, then the *cylinder* of $a$ in $T^\alpha$ is the set $\{r : \langle r_{k_1}, r_{k_2}, \ldots r_{k_n} \rangle \in a,\ r \in T^\alpha,\ \alpha \geq \omega\}$. Do not confuse the two concepts: cylinder of $a$ and the $y_k$-cylinder $C_{y_k} a$ of $a$.

The knowledge of the concept of the finite power of a $\sigma$-set algebra (basic set algebra) is assumed. The basic set algebra will be here the one-dimensional Borel-, Lebesgue $\sigma$-algebra, or a one-dimensional discrete $\sigma$-algebra.

If the dimension sets are finite for all the members of a set algebra $\mathcal{F}$ with unit $T^\alpha$, $\alpha \geq \omega$, then $\mathcal{F}$ is called *locally finite-dimensional* (in short, *locally finite*) set algebra. Let $\mathcal{B}^{cyl}$ and $\mathcal{L}^{cyl}$ denote the set algebras of the cylinders of the finite-dimensional Borel- or Lebesgue sets in $[0,1)^\alpha$, respectively (see [2]). These algebras are locally finite-dimensional algebras and not $\sigma$-algebras.

The operation $C^\partial_{y_k}$ can be defined in terms of $C_{y_k}$ by $C^\partial_{y_k} a = - C_{y_k}(-a)$. From the viewpoint of logic, the *universal quantification* $\forall y_k$ can be associated with $C^\partial_{y_k}$.

The transformation in which a free variable $y_m$ in $a \subset T^\omega$ is substituted by a variable $y_n$ is denoted by $S^{y_m}_{y_n} a$. $S^{y_m}_{y_n} a$ is also referred to as changing the dimension $y_m$ for $y_n$ in $a$. $S^y_z a$ is denoted also by $a(y/z)$. The $\alpha$-powers of set algebras are closed under changing dimensions. If a set algebra $\mathcal{F}$ with unit $T^\alpha, \alpha \geq \omega$ is closed under the operations cylindrifications $C_{y_k}$'s, $k < \alpha$ and under the substitutions $S^{y_m}_{y_n}$ then $\mathcal{F}$ is called an infinite-dimensional *polyadic set algebra* (see [4], [9], [5]).

The equality (1) is remarkable for our results. Assume that $\mathcal{D}$ is a locally finite-dimensional polyadic set algebra with unit $T^\alpha, \alpha \geq \omega$.

*Proposition* (1) *is valid under the conditions at* (1).
(see [10] I.1.11)

Proof. The definition of $C_y a$ implies that $a(y/y_i) \subset C_y a$ for every $i \in \omega$, i.e. $C_y a$ is an upper bound of the sets $a(y/y_i)$, $i \in \omega$. To prove that $C_y a$ is the least upper bound, indirectly, assume that a set $b \in D$ is the least upper bound, i.e. $b$ is an upper bound and $b \subset C_y a$ strictly. It is enough to prove that $C_y a = \sup_{i \in \omega} \{a(y/y_i)\}$, because $\sup_{i \in \omega} \{a(y/y_i)\} = \sup_{i \in \omega} \bigcup_{i=1}^n a(y/y_i)$. By condition, $a(y/y_i) \subset b$ for each $i \in \omega$. $D$ is locally finite, therefore there is a $y_k$ in the sequence $\langle y_i : i \in \omega \rangle$ such that $y_k \notin \triangle b \cup \triangle a$. For this $y_k$, $a(y/y_k) \subset b$ is true as well. Applying $C_{y_k}$ for this inclusion, $C_{y_k} a(y/y_k) \subset C_{y_k} b = b$ follows. But, $C_{y_k} a(y/y_k) = C_y a$ because $y_k \notin \triangle a$. Thus $C_y a \subset b$ follows and this contradicts $b \subset C_y a$ strictly, therefore $b = C_y a$.

Notice that (1) is valid even if $\mathcal{D}$ is not closed under cylindrifications but, specially, $C_y a \in \mathcal{D}$.

Considering (1), notice that if $a$ is finite-dimensional, then, in general, the infinite union $f = \bigcup_{i=1}^\infty a(y/y_i)$ does not belong to $D$ (where the union denotes here the union of sets). Because, on the one hand, by (1), $f = C_y a$, thus $f$ is a finite-dimensional element, on the other hand, $f$ as the infinite union $\bigcup_{i=1}^\infty a(y/y_i)$ is not finite-dimensional, in general. This is a contradiction. As a consequence, the sets $C_y a$ and $\bigcup_{i=1}^\infty a(y/y_i)$ are different.

There is an analog version of (1) for the operation $C^\partial_y a$, too. In the re-mainders, the supremum in (1) is called *cylindric sum*.

*2.2. Measures*

Let $\mathcal{F}$ be a set algebra with unit $T^\alpha$, where $T$ is an infinite set, $\alpha \geq \omega$. By a *measure* $\mu$ defined on $\mathcal{F}$ we mean a usual probaility measure, i.e.
(i) $0 \leq \mu(h) \leq 1$ for every $h \in F$
(ii) $\mu(\emptyset) = 0$, $\mu(T^\alpha) = 1$
(iii) $\sum_{i \in \omega} \mu(a_i) = \mu(a)$, if $\bigcup_{i \in \omega} a_i = a$, , $a_i \cap a_j = \emptyset$, $a, a_i \in F$.

Thus $\mu$ is $\sigma$-*additive* and $\mathcal{F}$ is not necessarily a $\sigma$-algebra.

The measure $\mu$ has the *product property* if $\triangle a \cap \triangle b = \emptyset$ implies that $\mu(a \cap b) = \mu(a) \cdot \mu(b)$. The measure $\mu$ is *symmetrical* if $\mu(a) = \mu(b)$ when $b$ is composed from $a$ by a substitution (by changing dimensions).

The knowledge of the concept of the finite power of a measure defined on a one-dimensional $\sigma$-set algebra is assumed ([2], [8]). It is known that the power measures have the product and symmetry properties.

A measure $\mu$ defined on a locally finite-dimensional polyadic set algebra $\mathcal{D}$ is *continuous* w.r.t. the cylindric sums in (1) if (2), i.e.

$$\mu(C_y a) = \sup_{n \in \omega} \mu(\bigcup_{i=1}^{n} a(y/y_i))$$

holds for every $a \in D$ and $y$, where the meaning of the notations is the same as in (1).

As it was mentioned, $C_y a \neq \bigcup_{i=1}^{\infty} a(y/y_i)$, thus (2) can not be reduced directly to some usual continuity of the measures.

The property analogous to (2) for $C_y^\partial a$ is:

$$\mu(C_y^\partial a) = \inf_{n \in \omega} \mu(\bigcap_{i=1}^{n} a(y/y_i)) \qquad (5)$$

The suprema and infima in (2) and (5) are called *cylindric suprema and cylindric infima*. We can speak about the continuity of the measure with respect to *only a given* supremum in (1).

*2.3. A basic lemma*

Let $\mathcal{D}$ be a locally finite-dimensional cylindric set algebra, $a \in D, C_y^\partial a \in D$ and let $\mu$ be a measure on $\mathcal{D}$.

**Lemma 1.** *Both the assumptions a) and b) below are necessary and sufficient for the validity of the equality in (5)*
  a)
$$\inf_{n\in\omega} \mu(\bigcap_{i=1}^n b\ (y/y_i)) = 0 \qquad (6)$$

*holds for arbitrary $\omega$-sequence of distinct variables $y_1, y_2, \ldots y_i, \ldots$ not in $\triangle a$, $y \in \triangle a$, where $b$ denotes the set $a - C_y^\partial a$.*
  b) $C_y^\partial e = \emptyset, e \in D$ *implies that*

$$\inf_{n\in\omega} \mu(\bigcap_{i=1}^n e(y/y_i)) = 0 \qquad (7)$$

*for arbitrary $\omega$-sequence of distinct variables $y_1, y_2, \ldots y_i, \ldots$ not in $\triangle e$, $y \in \triangle e$.*

**Proof.** a) The necessity is obvious because $C_y^\partial(a - C_y^\partial a) = C_y^\partial(a \cap -C_y^\partial a) = \emptyset$. By (5), $0 = \mu(\emptyset) = \inf_{n\in\omega} \mu(\bigcap_{i=1}^n b(y/y_i))$.

The sufficiency:
$(a - C_y^\partial a)(y/y_i) = (a \cap -C_y^\partial a)(y/y_i) = a(y/y_i) \cap -C_y^\partial a$ thus, $\bigcap_{i=1}^n (a \cap -C_y^\partial a)(y/y_i) = \bigcap_{i=1}^n a(y/y_i) \cap -C_y^\partial a$. By (6),

$$\inf_{n\in\omega} \mu[\bigcap_{i=1}^n a(y/y_i) \cap -C_y^\partial a] = 0 \qquad (8)$$

By definition, $C_y^\partial a \subset \bigcap_{i=1}^n a(y/y_i)$, therefore $\mu(\bigcap_{i=1}^n a(y/y_i)) =$
$= \mu(C_y^\partial a) + \mu(\bigcap_{i=1}^n a(y/y_i) \cap -C_y^\partial a)$. Applying the $\inf_{n\in\omega}$ to the latter sum $\inf_{n\in\omega} \mu(\bigcap_{i=1}^n a(y/y_i)) = \mu(C_y^\partial a) + \inf_{n\in\omega} \mu(\bigcap_{i=1}^n a(y/y_i) \cap -C_y^\partial a)$, using (8), we get that $\inf_{n\in\omega} \mu(\bigcap_{i=1}^n a(y/y_i) = \mu(C_y^\partial a))$ really, as we have to prove.

b) The necessity is obvious. The sufficiency holds because the condition implies (6), choosing $a - C_y^\partial a$ for $e$. ∎

Regarding this Section see also the references [6], [7], [3].

## 3. Discrete measures

Let $T$ be a fixed countable subset of $R^1$, i.e. $T = \{c_i\}_{i \in \omega} \subset R^1$ and $\alpha \geq \omega$. Let $p$ be a fixed, discrete measure on $T$. Let us consider $\mathcal{P}(T)$ as a $\sigma$-set algebra. For fixed $n \in \omega$ let us consider the $n$th power $\sigma$-set algebra $(\mathcal{P}(T))^n$. The $n$th power $p^n$ is a $\sigma$-additive measure on $(\mathcal{P}(T))^n$. If an infinite ordinal $\alpha$ is fixed, then the cylinders of the elements of the possible finite powers of $\mathcal{P}(T)$ form a set algebra $\mathcal{H}$ with unit $T^\alpha$ and the possible finite powers of $p$ define a $\sigma$-additive measure $v$ on this algebra ([2]).

**Theorem 2.** $\mathcal{H}$ equipped with the cylindrifications and with the substitutions $y/y_i$ in $T^\alpha$ is a polyadic set algebra. The measure $v$ defined on $\mathcal{H}$ is continuous, i.e. it satisfies the property (2).

**Proof.** The cylindrification of an $a \in \mathcal{H}$ is also an element of $\mathcal{H}$ because the definition of the finite powers of $\sigma$-algebras. So $\mathcal{H}$ is closed under cylindrifications. $\mathcal{H}$ is obviously closed under the substitutions $y/y_i$ because $\mathcal{H}$ is defined by power set algebras. Thus $\mathcal{H}$ can be regardeded as a polyadic set algebra.

$\mathcal{H}$ is locally finite, thus (1) holds in $\mathcal{H}$. We prove (5) rather than (2), and (7) rather than (5), using Lemma 1. Let $e$ be $a$.

For simplicity, let us suppose that the size of the dimension set of $a$ is 2, say $\triangle a = \{y, z\}$. Assume that $C_y^\partial a = \emptyset$.

If $n = 1$ in (7), then $v(a(y,z)) = \sum_{c_i \in T} v(a(y, c_i))$ by the $\sigma$-additivity of $v$.

If $n = 2$ in (7), then
$v(a(y,z) \cap a(x,z)) = \sum_{c_i \in T} v(a(y, c_i) \cap a(x, c_i))$ by the $\sigma$-additivity of $v$, where $y_1 = y$ and $y_2 = x$ in (7). However, $\triangle a(y, c_i) \cap \triangle a(x, c_i) = \emptyset$, so by the product property of $v$, $v(a(y, c_i) \cap a(x, c_i)) = v(a(y, c_i)) \cdot v(a(x, c_i))$.

The symmetry property of $v$ implies that the members of the product are equal, therefore $v(a(y,z) \cap a(x,z)) = \sum_{c_i \in T} v^2(a(y, c_i))$.

Generalizing this calculation from $n = 2$ to $n$ and to the sequence $y_1, y_2, \ldots y_n$ we get that in (7)

$$v(\bigcap_{i=1}^{n} a(y/y_i)) = \sum_{c_i \in T} v^n(a(y, c_i)) \quad (9)$$

The assumption $C_y^\partial a(y, z) = \emptyset$ implies that $v(a(y, c_i)) < 1$ because $v$ is a probability measure. Therefore, for fixed $c_i$, $v^n(a(y, c_i)) \to 0$ if $n \to \infty$.

We check that the infimum in question in (7) is 0, i.e. the infimum in the summa in (9) is 0. It is proven that for every $\varepsilon > 0$ there exists $M$, such that if $n > M$, then the summa in (7) less than $\varepsilon$.

The summa in (9) exists for any fixed $n$. Hence we can choose such an $N$, that $\sum_{i=N+1}^{\infty} v^n(a(y, c_i)) < \frac{\varepsilon}{2}$. Here the value of the left hand side is decreasing if $n$ increases. For fixed $N$, $\sum_{i=1}^{N} v^n(a(y, c_i)) \to 0$ if $n \to \infty$, therefore, for $\frac{\varepsilon}{2}$ and suitable $M$ if $n > M$, then $\sum_{i=1}^{N} v^n(a(y, c_i)) < \frac{\varepsilon}{2}$. Therefore $\sum_{i=1}^{N} v^n(a(y, c_i)) + \sum_{i=N+1}^{\infty} v^n(a(y, c_i)) < \frac{\varepsilon}{2} + \frac{\varepsilon}{2} = \varepsilon$ if $n > M$.

If $a$ has $k$ free variables besides $z$, say the variables $z_2, z_3, \ldots z_k$, i.e. $\triangle a = \{y, z, z_2, z_3, \ldots z_k\}$, then introducing the vector $\overline{z} = \{z, z_2, z_3, \ldots z_k\}$, the above proof can be repeated for the set $a(y, \overline{z})$. ∎

The $\sigma$-additivity of $\nu$ implies that $v(\bigcup_{i=1}^{\infty} a(y/y_i)) = \sup_{n \in \omega} v(\bigcup_{i=1}^{n} a(y/y_i))$. But, the above theorem implies that $v(C_y a) = \sup_{n \in \omega} v(\bigcup_{i=1}^{n} a(y/y_i))$. Consequently $v(C_y a) = v(\bigcup_{i=1}^{\infty} a(y/y_i))$ though $C_y a \neq \bigcup_{i=1}^{\infty} a(y/y_i)$.

## 4. Lebesgue measure

Let us consider the $\alpha$-power $\langle [0, 1)^\alpha, (\mathcal{L}_1)^\alpha, (\lambda_1)^\alpha \rangle$ $(\alpha \geq \omega)$ of the one-dimensional Lebesgue measure space $\langle [0, 1), \mathcal{L}_1, \lambda_1 \rangle$. It is denoted also by $\langle [0, 1)^\alpha, \mathcal{L}, \lambda \rangle$, where $\lambda$ is the Lebesgue measure on $\mathcal{L}$.

It is known that the finite unions of the cylinders of the finite-dimensional rectangles in $[0, 1)^\alpha$ form a set algebra $\mathcal{G}$, where the sides of these rectangles are non-degenerate subintervals of $[0, 1)$ ([8], [2]).

**Theorem 3.** *The set algebra $\mathcal{G}$ with the operations cylindrifications and the substitutions $y/y_i$ in $T^\alpha$ forms a polyadic set algebra with unit $[0,1)^\alpha$. The Lebesgue measure $\lambda$ on $\mathcal{G}$ is continuous, i.e. it satisfies the property (2).*

**Proof.** $\mathcal{G}$ is closed under any cylindrifications $C_y$ because if an element of $a$ is of the form $\bigcup_{j=1}^{k} I_j$, where the $I_j$'s are finite-dimensional rectangles, then $C_y \bigcup_{j=1}^{k} I_j = \bigcup_{j=1}^{k} C_y I_j$ and $C_y I_j$ is again a finite-dimensional rectangle. $\mathcal{G}$ is closed under changing dimensions because the definition of the rectangle. Thus, $\mathcal{G}$ can be regarded as a locally finite polyadic set algebra $D$. We prove the continuity of $\lambda$ w.r.t. the cylindric sums. We are going to use the infimum versions (5) and (7) in Lemma 1. Let $e$ be $a$.

First, assume that $\alpha$ has only one additional free variable $z$ other than $y$, so $\alpha$ is of the form $\alpha(y, z)$, i.e. $a$ is 2-dimensional. Later, this case can be generalised for finitely many free variables. It is known that the finite unions of finite-dimensional rectangles can be considered as finite *disjoint* unions of finite-dimensional rectangles. Let us consider the set $a(y, z)$ in this form.

Assume that $n = 1$ in (7).

Let us take the projection of $a(y, z)$ on the axis $z$. This is a finite union of one-dimensional rectangles. This union can be regarded as a disjoint union of some rectangles $I_1, I_2, \ldots I_m$. Then the union of the original 2-dimensional disjoint union can be composed in the form $\bigcup_{j=1}^{m} \bigcup_{i=1}^{i(j)} (I_j \times K_i^j)$, where $K_i^j$'s are one-dimensional disjoint intervals. Then

$$\lambda(\alpha(y,z)) = \lambda\left(\bigcup_{j=1}^{m} \bigcup_{i=1}^{i(j)} (I_j \times K_i^j)\right) = \sum_{j=1}^{m} \sum_{i=1}^{i(j)} \lambda(I_j) \times \lambda(K_i^j) =$$

$$= \sum_{j=1}^{m} \lambda(I_j) \cdot \lambda(\bigcup_{i=1}^{i(j)} K_i^j),$$ where $\lambda$ is the Lebesgue measure. Denote $q_j$ the value $\lambda(\bigcup_{i=1}^{i(j)} K_i^j)$.

By the assumption (7), $C_y^\partial \alpha(y, z) = \emptyset$. Therefore, for any $j$ and $z \in I_j$, there is a point $(y_0, z)$ not included in $\alpha(y, z)$ thus, $\lambda(\bigcup_{i=1}^{i(j)} K_i^j) < 1$, i.e. $q_j < 1$.

If $n = 2$ in (7), let us take the set $\alpha(y_1, z) \cap \alpha(y_2, z)$. In the space $(y_1, y_2, z)$

$$\lambda(\alpha(y_1,z) \cap \alpha(y_2,z)) = \lambda\left(\bigcup_{j=1}^{m}\bigcup_{i=1}^{i(j)}\bigcup_{p=1}^{p(j)}(I_j \times K_i^j \times H_p^j)\right) =$$

$$= \sum_{j=1}^{m}\sum_{i=1}^{i(j)}\sum_{p=1}^{p(j)} \lambda(I_j \times K_i^j \times H_p^j) = \sum_{j=1}^{m} \lambda(I_j) \cdot \lambda(\bigcup_{i=1}^{i(j)} K_i^j) \cdot \lambda(\bigcup_{p=1}^{p(j)} H_p^j) =$$

$$= \sum_{j=1}^{m} \lambda(I_j) \cdot q_j^2, \text{ where } q_j < 1, \; i(j) = p(j) \text{ and } \lambda(\bigcup_{i=1}^{i(j)} K_i^j) = \lambda(\bigcup_{p=1}^{p(j)} H_p^j), \text{ as a}$$

result of $C_y^\partial \alpha(y,z) = \emptyset$, the product property and the symmetry of $\lambda$.

More generally,
$$\lambda(\bigcap_{i=1}^{n} \alpha(y_i,z)) = \sum_{j=1}^{m} \lambda(I_j) \cdot q_j^n$$

where $q_j < 1$. This implies that the limit of the right-hand side is 0 for $n \to \infty$, hence the infimum of $\lambda(\bigcap_{i=1}^{n} \alpha(y_i,z))$ $n \in \omega$, is also 0, as we have to prove.

If $a$ has two additional variables $z_1$ and $z_2$ besides $y$, i.e. $a$ is 3-dimensional, the proof of (2) is similar. The first step is to take the projection of $a$ on the plane $(z_1, z_2)$, then we take the disjoint union composition of the projection. $a$ can be composed in the form $\bigcup_{j=1}^{m}\bigcup_{i=1}^{i(j)}(R_j \times K_i^j)$ for some $m$ and $i(j)$, where $R_j$'s are disjoint rectangles on the plane $(z_1, z_2)$ and $K_i^j$'s are one-dimensional disjoint intervals. The sequel of the proof is the same as above.

The proof of the $n$-dimensional case is analogous to those of the above cases. ∎

1. Another possible proof for Theorem 3 is to reduce the theorem to Theorem 2. The idea here that $\lambda$ is a strictly positive measure on the Boolean set algebra $\mathcal{G}$ satisfying the so called *countable-chain condition* in Boolean algebras.

2. Let us consider the cylinders of the so called *Borel rectangles* $B_1 \times B_2 \times \ldots \times B_n$ in $[0,1)^\alpha$, where $n$ runs over $\omega$, and $B_1, B_2, \ldots B_n$ are arbitrary one-dimensional Borel sets.

A generalization of the first part of Theorem 3 is that the possible finite unions of the cylinders of the Borel rectangles form a polyadic set algebra too, because $C_k(\bigcup_{i=1}^{n} F_i) = \bigcup_{i=1}^{n} C_k F_i$ holds, where the $F_i$'s are Borel rectangles, and $C_k F_i = B_k$ for some $k$, if $B_k$ is included in $F_i$, else $C_k F_i = [0,1)$. However, the second part of Theorem 3 is false for Borel rectangles as the example below shows.

Let $a(y,z)$ be the union of $F_1$ and $F_2$, where $F_1 = \{\frac{2}{3}\} \times [\frac{1}{2}, 1)$ and $F_2 = \{\frac{1}{3}\} \times [0, \frac{1}{2})$. Then $C_y a(y,z) = [0,1)$, $\lambda(C_y a) = 1$ and $\lambda(F_1 \cup F_2) = \lambda(a(y,z)) = 0$. By the symmetry of $\lambda$, on the right-hand side $\lambda(\bigcup_{i=1}^{n} a(y_i, z))$ of (2) is also 0. Thus (2) is false for $\lambda$.

In general, if $a$ is a Borel set, the continuity of $\lambda$ w.r.t. cylindric sums or products, not necessarily true. An example is: let $a(y,z) = \{\langle y, z\rangle : y = z\}$. Then $C_y a(y,z) = [0,1)$, but $\lambda(a(y,z) = 0$. $\lambda(C_y a) = 1$. By the symmetry of $\lambda$, $\lambda(\bigcup_{i=1}^{n} a(y_i, z)) = 0$.

Next a sufficient condition is given for the continuity of the Lebesgue measure $\lambda$ w.r.t. cylindric sums. The idea behind it is of geometrical character: $a$ can be included into the majorant rectangle set $C_{z_1} C_{z_2} \ldots C_{z_m} a$.

Assume that $a \in \mathcal{B}^{cyl}$, $\triangle a = \{z_1, z_2, \ldots z_m, y\}$ and the sets $a$, $C_{z_2} a$, $C_{z_1} C_{z_2} a, \ldots C_{z_1} C_{z_2} \ldots C_{z_n} a$ are element of $\mathcal{L}^{cyl}$.

**Theorem 4.** $\lambda(C_{z_1} C_{z_2} \ldots C_{z_m} a) < 1$ implies the sup property in (2).

**Proof.** First, let us assume that $\triangle a = \{y, z\}$, so besides $y$, $a$ has only one free variable $z$.

The product property applies to the product $\bigcap_{i=1}^{n} C_z a(y/y_i)$, thus

$\lambda(\bigcap_{i=1}^{n} C_z a(y/y_i)) = \prod_{i=1}^{n} \lambda(C_z a(y/y_i))$.

But, by the symmetry property, $\lambda(a(y/y_i) = \lambda(a(y/y_j), i \neq j$. Thus $\prod_{i=1}^{n} \lambda(C_z a(y/y_i)) = [\lambda(C_z a(y,z))]^n$.

It is shown that $\inf_{n} \lambda(\bigcap_{i=1}^{n}(a(y/y_i))) = 0$. $\bigcap_{i=1}^{n} a(y/y_i) \subset \bigcap_{i=1}^{n} C_z a(y/y_i)$, so $\lambda(\bigcap_{i=1}^{n} a(y/y_i)) \leq \lambda(\bigcap_{i=1}^{n} C_z a(y/y_i))$. But, $\lambda(\bigcap_{i=1}^{n} C_z a(y/y_i)) = [\lambda(C_z a(y,z))]^n$. Thus $\lambda(C_z a) < 1$ implies that $\inf_{n} \lambda(\bigcap_{i=1}^{n}(a(y/y_i))) = 0$, really. This follows also in that case too, when especially, $C_y^\partial a = \emptyset$. We can use the argument in the Lemma. Therefore, (7) implies that (5) holds.

If $a$ has exactly $m$ free variables $z_1, z_2, \ldots z_m$ besides $y$, then we can repeat the argument above $m$ times. ∎

Notice that the condition in the theorem is only sufficient, see the case when $\lambda(a) = 1$.

## 5. Strong projections

As is known, the projections of a finite-dimensional Borel- or Lebesgue measurable set is not necessarily Borel- or Lebesgue measurable, repectively. Therefore, $\mathcal{B}^{cyl}$ or $\mathcal{L}^{cyl}$ cannot be regarded as a polyadic set algebra. The projection of a Borel set is Lebesgue measurable (as a consequence of the *measurable projection theorem*, see [12], [3]), but, as the examples show in Sect. 4, the continuity of the Lebesgue measure is not necessarily true for cylindric sums in the sense of (2).

The concept of *strong projection* is introduced. $\mathcal{L}^{cyl}$ will be closed under this type of projection, and the Lebesgue measure $\lambda$ will be continuous w.r.t. these projections. Besides this, we present necessary and sufficient conditions for the Lebesgue measure to be continuous w.r.t. the usual projection and a given cylindric sum, i.e. for the satisfaction of (2).

Assume that the set algebra is $\mathcal{L}^{cyl}$, i.e. the set algebra of the cylinders of the finite-dimensional Lebesgue measurable sets, assume that $a \in \mathcal{L}^{cyl}$ and the measure is the Lebesgue measure $\lambda$ on $\mathcal{L}^{cyl}$.

First, assume that $\triangle a = \{y_k, z\}$, i.e. $a$ is 2-dimensional. The definition of the operation $\widetilde{C_{y_k}}$, corresponding to that of $C_{y_k}$, is as follows:

$$\widetilde{C_{y_k}} a(y_k, z) = \{r : \lambda(a(y_k, z) \mid z = r_z) > 0, r \in [0,1)^\alpha\} \tag{10}$$

and the definitionof the operation $\widetilde{C^\partial_{y_k}}$, corresponding to $C^\partial_{y_k}$, is as follows:

$$\widetilde{C^\partial_{y_k}} a(y_k, z) = \{r : \lambda(a(y_k, z) \mid z = r_z) = 1, r \in [0,1)^\alpha\} \tag{11}$$

where $\lambda(\mid z = r_z)$ denotes the one-dimensional Lebesgue conditional distribution w.r.t the condition $z = r_z$. The Lebesgue measure is power measure, so the conditional $\lambda(\mid z = r_z)$ is simply the one-dimensional Lebesgue measure, so $\lambda(a(y_k, z) \mid z = r_z) = \lambda(a(y_k, r_z))$, where.$a(y_k, r_z)$ is a one-dimensional set in the dimension $y_k$. In (10), $\lambda(a(y_k, z) \mid z = r_z) > 0$ is equivalent to $\lambda(a(y_k, z) \mid z = r_z) \neq 0$, of course.

$\widetilde{C_{y_k}}$ means forming a kind of $y_k$-cylinder set. The set $\widetilde{C_{y_k}}a$ is identified by the *strong projection* of $a$, parallel to $y_k$.

Like $C_{y_k}$ and $C^{\partial}_{y_k}$ the connection between $\widetilde{C_{y_k}}$ and $\widetilde{C^{\partial}_{y_k}}$ is

$$\widetilde{C_{y_k}}a = -(\widetilde{C^{\partial}_{y_k}} - a) \qquad (12)$$

If $\triangle a = \{y_k, z_{j_1}, \ldots z_{j_m}\}$, i.e. $a$ has $m$ free variables besides $y_k$, then the definition is completely similar to (10) and (11), but instead of $\lambda(a(y_k, z) \mid z = r_z)$,
$\lambda(a(y_k, z_{j_1}, \ldots z_{j_m}) \mid z_{j_1} = r_{j_1}, \ldots z_{j_m} = r_{j_m})$ is needed.
With the operations $\widetilde{C_{y_k}}$ and $\widetilde{C^{\partial}_{y_k}}$, certain special probability quantifiers can be associated in logic (see [13]). Notice that strong projection *depends* on the measure of the projected set.

Notice that in the case of discrete measures (Sect. 3) $\widetilde{C_{y_k}}$ and $C_{y_k}$ coincide. Therefore, the strong approach of assigning measure to a projection, in the case of Lebesgue measure can be regarded as a generalization of the discrete case.

**Lemma 5.** *The definitions of the operations $\widetilde{C_{y_k}}$ and $\widetilde{C^{\partial}_{y_k}}$ above make sense, i.e. if $a$ is Lebesgue measurable, then so are $\widetilde{C_{y_k}}a$ and $\widetilde{C^{\partial}_{y_k}}a$ thus, $\mathcal{L}^{cyl}$ is closed under the operations $\widetilde{C_{y_k}}$ and $\widetilde{C^{\partial}_{y_k}}$.*

**Proof.** Assume that $\triangle a = \{y_k, z\}$. The set $a$ is Lebesgue measurable, so its characteristic function is also Lebesgue measurable. The real function $\lambda(a(y_k, r_z), r_z)$ is Lebesgue measurable as the function $\lambda$ of $r_z$ (i.e. $\lambda(r_z)$) by Fubini's theorem (see [2]). This implies that the conditional probabilities in (10) and (11) exist because the measurability of $\lambda(r_z)$ implies that the sets $\{r_z : \lambda(r_z) > 0\}$ and $\{r_z : \lambda(r_z) = 1\}$ are Lebesgue measurable, so $\widetilde{C_{y_k}}a$ and $\widetilde{C^{\partial}_{y_k}}a$ exist and they are Lebesgue measurable. If $\triangle a$ is any finite set, the proof is completely analogous to the case above. ∎

The lemma implies that $\mathcal{L}^{cyl}$ with $\widetilde{C_{y_k}}a$ can be regarded as a Boolean algebra with operator.

**Theorem 6.** *Let $a \in \mathcal{B}^{cyl}$ and $\lambda$ be the Lebesgue measure on it. Then*

$$\lambda(\widetilde{C^{\partial}_y}\, a) = \inf_{n \in \omega} \lambda(\bigcap_{i=1}^{n} a(y/y_i)) \qquad (13)$$

*for arbitrary $\omega$-sequence of distinct variables $y_1, y_2, \ldots y_i, \ldots y_i \notin \triangle a$, $y \in \triangle a$.*

**Proof.** First, assume that $\triangle a = \{y, z\}$. By the previous lemma if $a \in \mathcal{B}^{cyl}$, then $\widetilde{C_y^\partial} a$ is Lebesgue measurable, so $\lambda(\widetilde{C_y^\partial} a)$ exists. Let $\lambda_n$ denote the $n$-dimensional Lebesgue measure.

$\lambda_2(a(y,z)) = \iint_N h(y,z) dy dz = \int_0^1 \lambda_2(a(y,z) \mid z) \, dz$ by Fubini's theorem ([2] 6.3), where $N = [0,1) \times [0,1)$, $h(y,z)$ is the characteristic function of $a(y,z)$, $\lambda_2(a(y,z) \mid z)$ denotes the Lebesgue conditional distribution w.r.t. the fixed condition $z$. But, $\lambda_2(a(y,z) \mid z) = \lambda_2(a(y,z))$ because the Lebesgue measure is a product measure.

Similarly, $\lambda_3(a(y_1,z) \cap a(y_2,z)) = \int_0^1 \lambda_2(a(y_1,z) \cap a(y_2,z) \mid z) dz =$

$= \int_0^1 \lambda_2(a(y_1,z) \cap a(y_2,z)) dz =$

$= \int_0^1 \lambda_1(a(y,z)) \cdot \lambda_1(a(y,z)) dz = \int_0^1 [\lambda_1(a(y,z))]^2 \, dz$, due to the product - and symmetry property of $\lambda$. In general

$$\lambda_{n+1}(\bigcap_{i=1}^n a(y/y_i)) = \int_0^1 [\lambda_1(a(y,z))]^n dz \qquad (14)$$

Let us consider the set $\widetilde{C_y^\partial} a(y,z)$ - denote it $A$ ($A \subset [0,1)$) - and decompose the integral in (14):

$\int_0^1 [\lambda_1(a(y,z))]^n dz = \int_A [\lambda_1(a(y,z))]^n dz + \int_{\sim A} [\lambda_1(a(y,z))]^n dz$

In the first member, by the definition of $\widetilde{C_y^\partial} a(y,z)$, the integrand $[\lambda_1(a(y,z))]^n$ is equal to 1, while in the second member $[\lambda_1(a(y,z))]^n < 1$ by definition )(see 11) . Thus,

$\lim_{n \to \infty} \int_0^1 [\lambda_1(a(y,z))]^n dz = \lim_{n \to \infty} (\int_A [\lambda_1(a(y,z))]^n dz + \int_{\sim A} [\lambda_1(a(y,z))]^n dz) =$

$$= \lambda_1(A) + \lim_{n \to \infty} \int_{\sim A} [\lambda_1(a(y,z) \mid z)]^n dz \qquad (15)$$

However, $\lim_{n \to \infty} \int_{\sim A} [\lambda_1(a(y,z))]^n dz = \int_{\sim A} \lim_{n \to \infty} [\lambda_1(a(y,z))]^n dz = 0$ because of the boundedness of the integrand.

Using the general notation $\lambda$ for the Lebesgue measure we get from (11), (14) and (15) that $\inf_{n\in\omega} \lambda \bigcap_{i=1}^{n} a(y/y_i) = \lambda(\widetilde{C_y^\partial}\, a)$.

If $a$ has more free variables besides $y$, the proof is similar. ∎

The equality concerning $\widetilde{C_y}a$ corresponding to (13) is the formula

$$\lambda(\widetilde{C_y}\, a) = \sup_{n\in\omega} \lambda(\bigcup_{i=1}^{n} a(y/y_i)) \tag{16}$$

by (12).

An immediate consequence of (16) is:

**Corollary 7.** *Assume that $a$ and $C_y a$ are in $\mathcal{L}^{cyl}$. A necessary and sufficient condition for the continuity of the Lebesgue measure $\lambda$ w.r.t. the cylindric sums in (2) is*

$$\lambda(C_y a) = \lambda(\widetilde{C_y}\, a) \tag{17}$$

A sufficient conditions for (17) is $\lambda a = 1$, or $\lambda(C_y a) = 0$, for example.

Taking into consideration (16) and the fact that in discrete case the usual and the strong projections coincide, and in the discrete case (2) is true (see Theorem 2), we can conlude that the *quantity*

$$\sup_{n\in\omega} \lambda(\bigcup_{i=1}^{n} a(y/y_i))$$

is the measure of the strong projection both in the case of the discrete measures and the Lebesgue measure.

The results in this section can be generalised to any continuous power measure space rather than the Lebesgue measure space.